\documentclass[12pt]{article}

\makeatletter
\newfont{\footsc}{cmcsc10 at 8truept}
\newfont{\footbf}{cmbx10 at 8truept}
\newfont{\footrm}{cmr10 at 10truept}
\makeatother
\usepackage{amsmath}
\usepackage{amssymb}
\usepackage{exscale}
\usepackage{amsthm}
\usepackage{epsfig}
\usepackage{verbatim}
\usepackage{setspace}

\theoremstyle{plain}
\newtheorem{theorem}{Theorem}[section]

\newtheorem{lemma}[theorem]{Lemma}
\newtheorem{corollary}[theorem]{Corollary}
\newtheorem{conjecture}[theorem]{Conjecture}

\theoremstyle{definition}
\newtheorem{definition}[theorem]{Definition}

\newtheorem{remark}[theorem]{Remark}

\DeclareMathOperator{\1}{\bf{1}}

\newcommand{\N}{{\mathbb{N}}}

\newcommand{\R}{{\mathbb{R}}}

\newcommand{\PP}{{\mathbb{P}}}

\newcommand{\E}{{\mathbb{E}}}
\newcommand{\Var}{{\mathbb{V}}}

\def\bs0{\bf 0}

\title{Discrepancy Bounds for a Class of Negatively Dependent Random Points Including Latin Hypercube Samples}

\author{Michael Gnewuch\thanks{Institut f\"ur Mathematik, Universit\"at Osnabr\"uck, 
Germany ({\tt michael.gnewuch@uni-osnabrueck.de}).}
\and Nils Hebbinghaus\thanks{Institut f\"ur Informatik, Christian-Albrechts-Universit\"at zu Kiel, Germany
({\tt nils.hebbinghaus@gmail.com}).}}

\begin{document}

\maketitle
\vskip 1pc

\begin{abstract}
We introduce a class of $\gamma$-negatively dependent random samples. We prove that this class includes, apart from Monte Carlo samples, in particular Latin hypercube samples and Latin hypercube samples padded by Monte Carlo.

For a  $\gamma$-negatively dependent $N$-point sample in dimension $d$ we provide probabilistic upper bounds for its star discrepancy with explicitly stated dependence on $N$, $d$, and $\gamma$. These bounds generalize the probabilistic bounds for Monte Carlo samples from [Heinrich et al., Acta Arith. 96 (2001), 279--302] and [C.~Aistleitner, J.~Complexity 27 (2011), 531--540], and they are optimal for Monte Carlo and 
Latin hypercube samples.
In the special case of Monte Carlo samples the constants that appear in our bounds improve substantially on the constants presented in the latter paper and in 
[C.~Aistleitner, M.~T.~Hofer, Math. Comp.~83 (2014), 1373--1381].
\end{abstract}

\section{Introduction}

Discrepancy measures are well established and play an important role in fields like computer graphics, experimental design, empirical process theory, learning theory and machine learning, random number generation, optimization (in particular, stochastic programming), and numerical integration or  stochastic simulation, see, e.g., \cite{ASYM16, CM14, Cha00, DKS13, DP10, DGW14, DT97, FW93, Gla03, Lem09, LP14, Mat10, Nie92, NW10} and the literature mentioned therein.

The prevalent and most intriguing discrepancy measure is arguably the \emph{star discrepancy}, which is defined in the following way:

Let $P\subset [0,1)^d$ be an $N$-point set. (We always understand an ``$N$-point set'' as a ``multi-set'', i.e., it consists of $N$ points, but those points do not have to be pairwise different.) We define the \emph{local discrepancy} of $P$ with respect to a Lebesgue-measurable test set $T\subseteq [0,1)^d$ by
\begin{equation*}
D_N(P,T) := \bigg| \frac{1}{N} |P\cap T| - \lambda^d(T) \bigg|,
\end{equation*}
where $|P\cap T|$ denotes the cardinality of the finite set $P\cap T$ and $\lambda^d$ denotes the $d$-dimensional Lebesgue measure on $\R^d$.
For  vectors $x= (x_1,x_2,\ldots, x_d)$, $y= (y_1, y_2, \ldots, y_d) \in \R^d$ 
we write 
\begin{equation*}
[x,y) := \prod^d_{j=1} [x_j,y_j) = \{z\in \R^d \,|\, x_j \le z_j < y_j \hspace{1ex}\text{for $j = 1,\ldots, d$}\}.
\end{equation*}
The \emph{star discrepancy} of $P$ is then given by
\begin{equation*}
D^*_N(P) := \sup_{y\in [0,1]^d} D_N(P,[0,y)).
\end{equation*}

The star discrepancy is intimately related to quasi-Monte Carlo integration via the Koksma-Hlawka inequality (\cite{Hla61, Kok42}): For every $N$-point set $P\subset [0,1)^d$ we have
\begin{equation*}
\left| \int_{[0,1)^d} f(x) \,{\rm d}\lambda^d(x) - \frac{1}{N} \sum_{p\in P} f(p) \right| 
\le D_N^*(P) {\rm Var}_{\rm HK}(f), 
\end{equation*}
where ${\rm Var}_{\rm HK}(f)$ denotes the variation of the integrand $f$ in the sense of Hardy and Krause, see, e.g., \cite{AD15, Nie92}. 
The Koksma-Hlawka inequality is sharp, see again \cite{Nie92}.
(An alternative sharp version of the Koksma-Hlawka inequality can be found in \cite{HSW04}; it says that the worst-case error 
of equal-weight quadratures based on a set of sample points $P$ over the norm unit ball of the Sobolev space of dominated mixed smoothness of order one is exactly the star discrepancy of $P$.)
The Koksma-Hlawka inequality shows that equal-weight quadratures based on sample points with small star discrepancy yield  small  integration errors. (Deterministic equal-weight quadratures are commonly called \emph{quasi-Monte Carlo quadratures}; for a survey we refer to \cite{DKS13}.)  
For the very important task of high-dimensional integration, which occurs, e.g., in computational finance, physics or quantum chemistry, it is therefore of interest to know tight bounds for the smallest achievable star discrepancy
\begin{equation*}
D^*(N,d) := \inf\{ D^*_N(P) \,|\, P \subset [0,1)^d, \, |P| = N\}
\end{equation*}
or, equivalently, for the inverse of the star discrepancy
$$
N^*(\varepsilon, d) := \inf \{ N\in \N_0 \,|\, D^*(N,d)\le \varepsilon \}, 
$$
which is the minimum number of sample points that guarantee a discrepancy bound of at most $\varepsilon$,
and to be able to construct integration points that satisfy those bounds. To avoid the ``curse of dimensionality'' it is crucial that such bounds scale well with respect to the dimension $d$.

For fixed $d$ the best known asymptotic upper bounds for $D^*(N,d)$ are of the form 

\begin{equation}\label{disc_asymptotic}
D^*(N,d) \le C_d \ln(N)^{d-1} N^{-1}, \hspace{3ex}N\ge 2,
\end{equation}
see \cite{Hal60} or, e.g., the books \cite{DP10, Nie92}.
For larger $d$ those bounds give us no helpful information for moderate values of $N$, since the function $f(N):= \ln(N)^{d-1} N^{-1}$ is increasing for $N\le e^{d-1}$. Additionally,  for $d\ge 3$ values of $N$ much larger than $e^{d-1}$ are needed before $f(N)$ is below the common ``Monte Carlo rate'' $1/\sqrt{N}$ -- in dimension $d=10$, e.g.,  we need  $N$ to be larger than $1.295 \cdot 10^{34}$.
Moreover, the constant $C_d$ may grow unfavorably as $d$ gets larger. 
Actually, it is known for some $N$-point constructions $P$ that the constant $C_{d}'$ in the representation
\begin{equation*}
D^*(P) \le (C_d'\ln(N)^{d-1} + o(\ln(N)^{d-1})) N^{-1} 
\end{equation*}
of (1) tends to zero as $d$ approaches infinity, see, e.g., \cite{Nie92, NX96, Ata04} or \cite{FL09}.
But the behavior of the ``whole constant'' $C_d$ in \eqref{disc_asymptotic} is unfortunately not known.
That is why the asymptotic upper bound \eqref{disc_asymptotic} is not helpful for high-dimensional integration and we have to look for pre-asymptotic bounds, i.e., bounds that give us useful information already for a moderate number of points $N$.

The best known upper and lower bounds for the smallest achievable star discrepancy with explicitly given dependence on the number of sample points~$N$ as well as on the dimension~$d$ are of the following form: On the one hand, for all $d,N \in \N$ there exists
an $N$-point set $P\subset [0,1)^d$ satisfying
\begin{equation}\label{hnww01}
D^*_N(P) \le C \sqrt{\frac{d}{N}}
\end{equation}
for some constant $C>0$, implying
\begin{equation}\label{hnww02}
N^*(\varepsilon,d) \le \left\lceil C^2 d \varepsilon^{-2} \right\rceil
\end{equation}
for all $d\in \N$, $\varepsilon \in (0,1]$.
On the other hand, there exist constants $c, \varepsilon_0 >0$ such that 
\begin{equation}\label{hin03}
N^*(\varepsilon,d) \ge cd \varepsilon^{-1}
\end{equation}
for all $0<\varepsilon \le \varepsilon_0$, $d\in \N$, showing that 
for all $N$-point sets $P\subset [0,1)^d$ necessarily 
\begin{equation}\label{hin04}
D^*_N(P) \ge \min \left\{\varepsilon_0, c \frac{d}{N} \right\}.
\end{equation}
Notice that \eqref{hnww02} and \eqref{hin03} show that the inverse of the star discrepancy depends
essentially linearly on the dimension $d$ (in the sense that we have an upper bound for it that depends linearly on $d$ and also a lower bound that depends linearly on $d$).
The upper bounds \eqref{hnww01} and \eqref{hnww02} were proved 
in \cite{HNWW01} without providing an estimate for the constant $C$.
An estimate was given in \cite{Ais11}, who showed that $C\le 9.65$. In the course of this paper we will improve his estimate to $C \le 2.5287$ (see Corollary~\ref{Cor_Prob_Est}), consequently decreasing 
the minimum number of sample points needed to guarantee a given discrepancy bound $\varepsilon$ in arbitrary dimension $d$ 
by a factor of more than $14$, cf. \eqref{hnww02}. 
All the results mentioned so far are based on probabilistic arguments and do not provide an explicit
(deterministic) point construction that satisfies \eqref{hnww01}.
The lower bounds \eqref{hin03} and \eqref{hin04} were established in \cite{Hin04}.
Notice that there is a gap between the upper and lower bounds  (\ref{hnww01}) and (\ref{hin04}), and \eqref{hnww02} and \eqref{hin03}, respectively. 

Already 15 years ago
 Heinrich posed the following problems in \cite[Problem~1 \& 2]{Hei03}:
\begin{itemize}
\item[(P1)] For each $N,d \in \N$ give a (deterministic) construction of an $N$-point set $P \subset [0,1]^d$ satisfying 
\eqref{hnww01} for some positive constant $C$ not depending on $N$ or $d$.
\item[(P2)] Does any of the various known (deterministic) constructions of low discrepancy point sets
satisfy an estimate like (\ref{hnww01})? 
\item[(P3)] 
Determine the order of the smallest possible star discrepancy as a function of the number of points $N$ and the dimension $d$.
\item[(P4)] Determine 
\begin{equation*}
\begin{split}
&\overline{\alpha} := \\
&\sup \big\{ \alpha \,|\, \,\exists c,k\ge 0\, \forall N,d\in \N\, \exists P \subset [0,1]^d: |P|=N \wedge
D^*_N(P) \le c\tfrac{d^k}{N^{\alpha}} \big\}.
\end{split}
\end{equation*}
(The so-called \emph{exponent of tractability of the star discrepancy} $\tau$ is related to $\overline{\alpha}$ in the following way: 
$\overline{\alpha} = 1/\tau$; see, e.g., \cite{HNWW01}.)
\end{itemize}
Similar problems were stated in \cite{NW08, NW10} as Open Problems~6, 7, and 42.

Problems (P1) to (P4) turned out to be very difficult to solve.
It is, for instance, obvious that problem (P3) is a very hard problem, since it contains the so-called great open problem of discrepancy theory 
to find the precise order of the smallest possible star discrepancy in $N$ for fixed dimension $d\ge 3$. But also the other problems turned out to be very difficult 
to answer and have not been solved so far.  
For problem (P4) it is known that $1/2 \le \overline{\alpha} < 0.9373$, due to bound \eqref{hnww01} and a lower bound on the exponent of tractability of the star discrepancy from \cite{Mat98}. 
The following conjecture is due to Wo\'zniakowski, cf. \cite{Hei03} and \cite[Open Problem 7]{NW08}. 

\vspace{1ex}

\begin{conjecture} [Wo\'zniakowski] 
$\overline{\alpha} = 1/2$.
\end{conjecture}

\vspace{1ex}

If this conjecture is true, then, due to the essentially linear dependence of the inverse of the star discrepancy on the dimension $d$, the estimate \eqref{hnww01} is actually the best possible bound (apart from logarithmic factors) that is polynomial in $N^{-1}$ as well as in $d$. 

One reason for the difficulty of Problems (P1) and (P2) is that already the problem of  calculating the star discrepancy of an arbitrary $N$-point set is $NP$-hard, see \cite{GSW09}, and, in the language of parametrized complexity theory, $W[1]$-hard, see \cite{GKWW12}.

Heinrich also posed weaker versions of Problem (P1) and (P2),
cf. \cite{Hei03}. Those weaker versions were at least formally solved with the help of derandomized algorithms that generate deterministic $N$-point sets $P$ that satisfy 
\begin{equation}\label{wb1} 
D^*_N(P) \le C_1 \sqrt{\frac{d}{N}} \sqrt{\ln(1+N)},
\end{equation} 
see \cite{DGS05, DG08,DGW10}, or 
\begin{equation}\label{wb2} 
D^*_N(P) \le C_2 \sqrt{\frac{d^3}{N}} \sqrt{\ln(1+N/d)},
\end{equation} 
see \cite{DGKP08, DGW09}, for some small constants $C_1, C_2$. Those algorithms derandomize probabilistic experiments, in which random point sets $P$ satisfying \eqref{wb1} or \eqref{wb2} with high probability are generated. As numerical experiments in \cite{DGW09, DGW10} showed, those algorithms work well in dimensions up to $d=21$, but for much larger dimension their running times are prohibitive. (For a more extensive discussion, see also \cite{Gne11}.)

To get closer to a solution of the problems stated by Heinrich, we propose to study the following 
related randomized 
problems:
\begin{itemize}
\item[(R1)] What kind of randomized point constructions 
satisfy~(\ref{hnww01}) in expectation and/or with high probability? 
\item[(R2)] Are there randomized point constructions that
satisfy~(\ref{hnww01}) in expectation and/or with high probability and have more evenly distributed lower dimensional projections or satisfy better asymptotic discrepancy bounds than Monte Carlo points?
\item[(R3)]
Are there randomized point constructions that lead to a better estimate than \eqref{hnww01} or that can even be used to disprove Wo\'zniakowski's conjecture?
\end{itemize}
We believe that the problems (R1), (R2), and (R3) are important and interesting in their own rights. 
Moreover, an answer to question (R1) or (R2) would draw us closer to a solution of problem (P1) 
(since we may derandomize promising randomized point constructions) and of problem (P2) 
(since we get a hint, which known deterministic constructions are worth to be examined closer).
An affirmative answer to question (R3) may lead, due to the probabilistic method,
to progress in problem (P3) and  in problem (P4).

Let us explain this a little bit more.
As mentioned, the upper bound~\eqref{hnww01} was proved via probabilistic arguments. 
Indeed, Monte Carlo points, i.e., independent random points uniformly distributed in $[0,1)^d$, satisfy this bound with high probability (cf. also Corollary~\ref{Cor_Prob_Est}).

In~\cite{Doe13} it was shown that the star discrepancy of 
Monte Carlo point sets $X$ behaves like the right hand side in (\ref{hnww01}). More precisely,  there exists a constant $K>0$ such that the expected star discrepancy of $X$ is bounded from below by
\begin{equation}\label{benji1}
\E [D^*_N(X)] \ge K \sqrt{\frac{d}{N}}
\end{equation}
and additionally we have the probabilistic discrepancy bound
\begin{equation}\label{benji2}
\PP \! \left( D^*_N(X) < K\sqrt{\frac{d}{N}} \right) \le \exp(-\Omega(d)).
\end{equation}
The upper bound~\eqref{hnww01} is thus 
sharp for Monte Carlo points, showing that they cannot be employed to improve the upper discrepancy bound \eqref{hnww01} or to disprove Wo\'zniakowski's conjecture.

Clearly, there are other random point constructions that look more promising than simple Monte Carlo point sets as, e.g., Latin hypercube samples, see \cite{MBC79}, scrambled  $(t,m,s)$-nets, see \cite{Owe95, Owe97a,Owe97b}, or randomly shifted lattice rules, see \cite{SKJ02}.
In general, these random points will not be stochastically independent, which raises the next problem:
\begin{itemize}
\item[(R4)] How to analyze random point constructions whose single points are \emph{not}
stochastically independent? 
\end{itemize}
Obviously, this problem is not only of interest for analyzing the star discrepancy, but for stochastic simulation in general. 

In this paper we want to address the problems (R1), (R2), and (R4). We proceed as follows. First we introduce a negative dependence property of random point sets that is based on negative orthant dependence and that we call $\gamma$-negative dependence, see Definitions \ref{Def_Neg_gamma_Dep} and \ref{Def_neg_dep}. 
This property allows us to use large deviation bounds of Hoeffding- and Bernstein-type, see Section~\ref{EXP}. In Section~\ref{LHS} we prove that the class of $\gamma$-negatively dependent random points contains, in particular, Latin hypercube samples and Latin hypercube samples padded by Monte Carlo in arbitrary dimension $d$,  see Theorem~\ref{Thm-d-LHS}. With Theorem~\ref{theo} we provide a generalization of Theorem~\ref{Thm-d-LHS} that may be used to verify $\gamma$-negative dependence for random point sets different from (padded or unpadded) Latin hypercube samples. 
In Section~\ref{PDB} we show that for each $\gamma = e^{\rho d}$, $\rho \in \N_0$, all 
{$\gamma$-negative} dependent random point sets $P$ satisfy  \eqref{hnww01} with high probability and, surprisingly, the constant $C=C(\gamma)$ depends only mildly on $\gamma$, see Theorem \ref{main_theo}. In particular, our result generalizes the results
obtained in \cite{Ais11, AH14, HNWW01}, since the latter results can be seen as probabilistic discrepancy bounds for Monte Carlo point sets.   
In Corollary~\ref{Cor_Prob_Est} we provide  discrepancy bounds with explicit constants and explicit success probabilities for the special instances of Monte Carlo point sets and Latin hypercube samples (padded by Monte Carlo or not). In the special case of Monte Carlo point sets these constants and success probabilities improve substantially on the ones derived in \cite{Ais11} and \cite{AH14}.
In Remark~\ref{Rem_DDG18}
we explain that the probabilistic discrepancy bound \eqref{hnww01} is actually sharp for Latin hypercube samples; this result follows directly from new lower probabilistic discrepancy bounds in  \cite{DDG18}.

Let us close the introduction with some remarks on our notation.

For $N\in \N$ we denote the set $\{1,2,\ldots,N\}$ by $[N]$. We denote the Lebesgue measure on $\R$ by $\lambda$.  
By $\PP$, $\E$, and $\Var$ we always mean
probability, expectation, and variance, respectively.
If not specified otherwise, all random variables are defined on a 
probability space $(\Omega, \Sigma, \PP)$.

\section{Negative Dependence and Exponential Inequalities}
\label{EXP}

The concept of negative dependence was introduced in
\cite{Leh66} for pairs of random variables. 
In the literature one finds several contributions on rather demanding
notions of negative dependence as, e.g., negative association
introduced in \cite{JDP83}; a survey can be found in \cite{Pem00}.
Sufficient for our purpose is the following notion
for Bernoulli or binary random variables, i.e., 
random variables that only take
values in $\{0,1\}$.

\begin{definition}\label{Def_Neg_gamma_Dep}
Let $\gamma \ge 1$. 
We call binary random variables 
$T_1,T_2$, $\ldots$, $T_N$ 
\emph{upper $\gamma$-negatively dependent} if  
\begin{equation}
\label{cond2'}
\PP\left( \bigwedge_{j\in u} T_j = 1 \right) 
\le \gamma \prod_{j\in u} \PP( T_j = 1)
\hspace{2ex}\text{for all $u \subseteq [N]$,}
\end{equation}
and \emph{lower $\gamma$-negatively dependent} if  
\begin{equation}
\label{cond2''}
\PP\left( \bigwedge_{j\in u} T_j = 0 \right) 
\le \gamma \prod_{j\in u} \PP(T_j = 0)
\hspace{2ex}\text{for all $u \subseteq [N]$.}
\end{equation}
We call $T_1,T_2,\ldots,T_N$ \emph{$\gamma$-negatively dependent} if
both conditions (\ref{cond2'}) and (\ref{cond2''}) are satisfied.
If $\gamma = 1$, we usually will suppress the explicit reference to 
$\gamma$.
\end{definition}

A similar notion, called $\lambda$-correlation, can be found in 
\cite{PS97}.
 $1$-negative dependence is usually called negative orthant dependence, cf. \cite{BSS82}.
 
Notice that, in particular, independent binary random variables are negatively dependent. 
Furthermore, it is easily seen that for $N=2$ and $\gamma =1$ 
the notions of upper and lower $\gamma$-negative dependence are 
equivalent, cf. \cite{Leh66}.

We are interested in binary random variables $T_i$, $i=1, \ldots,N$, of the form
$T_i = \1_A(X_i)$, where $A$ is a Lebesgue-measurable subset of $[0,1]^d$ 
(whose characteristic function
is denoted by $\1_A$),  and $X_1, \ldots, X_N$ are randomly chosen points in $[0,1]^d$.
 
Panconesi and
Srinivasan derived in \cite{PS97} Chernoff-Hoeffding-type bounds (\cite{Che52, Hoe63})
for $\lambda$-correlated random variables. We will use the following two similar bounds of Hoeffding- and of Bernstein-type; for a proof see, e.g., \cite{Heb12}. 

\begin{theorem}\label{Hoeffding}
Let $\gamma \ge 1$, and let $T_1, \ldots , T_N$ be
$\gamma$--negatively dependent binary random variables. Put $S:= \sum^N_{i=1} (T_i-\E[T_i])$.
We have
\begin{eqnarray}\label{Hoeffding_both}
\mathbb{P}\left( |S| \geq t\right)  \leq 2\gamma\exp\left(-\frac{2t^2}{N}\right)
\hspace{2ex}\text{for all $t>0$.}
\end{eqnarray}
\end{theorem} 

\begin{theorem}\label{Bernstein}
Let $\gamma\ge 1$, and let $T_1,\ldots,T_N$ be $\gamma$-negatively dependent binary random variables. 
Put $S:= \sum^N_{i=1} (T_i-\E[T_i])$ and 
$\sigma^2 := \frac{1}{N} \sum^N_{i=1} \Var [T_i]$. 
Then we have
\begin{equation}
\label{bernstein}
\PP \left(|S| \ge t \right) \le 
2\gamma \exp \left( - \frac{t^2}{2N\sigma^2 + 2t/3} \right)
\hspace{2ex}\text{for all $t>0$.}
\end{equation}
\end{theorem}

 Let us close this section by mentioning a line of research that has been initiated in \cite{Lem18} and that is related to the one pursued in this paper. Its main goal is to show that for important classes of random variables the  variance of (suitable) randomized quasi-Monte Carlo estimators is never worse than the variance of the plain Monte Carlo estimator. The proof techniques used there are based on a pairwise negative dependence condition and Hoeffding's lemma, see \cite{Lem18}. Further results in this direction are provided in \cite{WL19, WG19, WGH19}.

\section{Latin Hypercube Sampling and Padding}
\label{LHS}

To provide useful examples of non-trivial $\gamma$-negatively dependent random variables, we prove in this section a 
negative dependence result for sample points stemming from a Latin hypercube sample, which may additionally  be padded by Monte Carlo. Further examples are provided in the follow-up paper \cite{WGH19}.
    
\begin{definition}
A \emph{Latin hypercube sample} (LHS) $(X_n)_{n=1}^N$ in $[0,1)^d$ is of the form
\begin{equation*}
X_{n,j} = \frac{\pi_j(n-1) + U_{n,j}}{N},
\end{equation*}
where $X_{n,j}$ denotes the $j$th coordinate of $X_n$, $\pi_j$ is a permutation
of $0,1$, $\ldots$, $N-1$, uniformly chosen at random, and $U_{n,j}$ is uniformly distributed in $[0,1)$.
The $d$ permutations $\pi_j$ and the $dN$ random variables $U_{n,j}$ are mutually independent. 
\end{definition}

The definition of Latin hypercube sampling presented above was introduced in  \cite{MBC79} for the design of computer experiments;
there is an older version of Latin hypercube sampling presented in \cite{Pat54}, where all uniform random variables $U_{n,j}$ are simply replaced by the value $0.5$.

As shown in \cite{Owe97b}, an estimator $\hat{I}$ of an integral $I$ based on Latin hypercube samples with $N$ points never leads to a variance greater than that of the corresponding  estimator based on $N-1$ Monte Carlo points.

\begin{remark}\label{1_dim_projections}
Notice that the one-dimensional projections of Latin hypercube samples are much more evenly distributed than the one-dimensional projections of Monte Carlo points.
This observation can be put into a quantitative statement by comparing the star discrepancies of the former and the latter projections: For a one-dimensional Latin hypercube sample of size $N$ the star discrepancy is at most $1/N$,
while for one-dimensional Monte Carlo samples of the same size the star discrepancy is 
of size $1/\sqrt{N}$, cf. \eqref{benji1} 
and \eqref{benji2}.
\end{remark}

\begin{definition}
Let $d,d',d''\in\N$ with $d=d'+d''$. Let $Y=(Y_k)_{k\in\N}$ be a 
(deterministic or randomized) sequence in $[0,1)^{d'}$, and
let $U=(U_k)_{k\in\N}$ be a sequence of independent uniformly 
distributed random vectors in $[0,1)^{d''}$. The $d$-dimensional concatenated 
sequence $X=(X_k)_{k\in \N} = (Y_k, U_k)_{k\in \N}$ 
is called a \emph{mixed sequence}. One also says that $X$ results from $Y$ by
\emph{padding by Monte Carlo}. 
\end{definition}

Padding by Monte Carlo was introduced in \cite{Spa95} to tackle
difficult problems in particle transport theory. He suggested to use a mixed sequence resulting from padding a deterministic low-discrepancy sequence.
Mixed sequences showed a favorable performance in several numerical experiments, see, e.g.,
\cite{Okt96, OTB06}. The latter papers also provided 
theoretical results on probabilistic
discrepancy estimates of mixed sequences which have been improved in
\cite{AH12, Gne09}.
Padding by LHS (instead of by Monte Carlo) was considered earlier in \cite[Example~5]{Owe94}.

We define the $1$-dimensional grid $G_N^1$ by
$$G_N^1:= \{ 0,1/N, \ldots, (N-1)/N, 1\}.$$  
The following lemma on $1$-dimensional LHS is the key ingredient in the proof of our main result on $d$-dimensional LHS (with or without padding by Monte Carlo), Theorem~\ref{Thm-d-LHS}. Theorem~\ref{Thm-d-LHS} combined with Theorem~\ref{main_theo} immediately imply the discrepancy bounds for LHS (with or without padding by Monte Carlo) in Corollary~\ref{Cor_Prob_Est}.

\begin{lemma}\label{one-d} 
Let $(X_n)_{n=1}^N$ be a LHS in $[0,1)$. Let $a,b\in  [0,1)$ with $a\le b$.
\begin{itemize} 
\item[{\rm (a)}] We have for all $\nu\in \{0,1, \ldots, N\}$ that 
\begin{equation}\label{pre*}
 \PP \left( \bigwedge_{\ell=1}^\nu X_{\ell} \in [a,b) \right) 
 \le (b-a)^{\nu}.
\end{equation}
\item[{\rm (b)}] Let $I^{(0)}_1 := [a,b)$, $I^{(0)}_2 := [0,b)$ and $I^{(1)}_1 := [b,1)$, $I^{(1)}_2 := [0,a) \cup [b,1)$.
Then we have for all $\sigma\in \{0,1\}$, $k\in \{0,1, \ldots, N\}$, and $\nu \in \{ 0,1, \ldots, k\}$ that
\begin{equation}\label{*}
 \PP \left( \bigwedge_{\ell=1}^\nu X_{\ell} \in I^{(\sigma)}_1 \wedge  \bigwedge_{\ell=\nu+1}^k X_{\ell} \in I^{(\sigma)}_2 \right) 
 \le \delta \lambda( I^{(\sigma)}_1)^{\nu} \lambda( I^{(\sigma)}_2)^{k-\nu}, 
 \end{equation}
 where 
 \begin{equation*}
 \delta := \begin{cases} 1 & \mbox{if} \quad a, b \in G^1_N \quad \mbox{or} \quad a=0,\\ 
 e & \mbox{else}.\end{cases}
\end{equation*}
The constant $\delta$ in (\ref{*}) is optimal in the following sense: for each $\delta < e$ and any $\sigma\in \{0,1\}$ there
exist $N\in \N$, $a,b\in [0,1)$, $k\in \{0,1, \ldots, N\}$, and $\nu \in \{ 0,1, \ldots, k\}$ such that (\ref{*}) does not hold.
\end{itemize}
\end{lemma}

\begin{proof}
Let $\alpha:= \lceil Na \rceil$, $\beta := \lfloor Nb \rfloor$, and let $\varepsilon_a, \varepsilon_b \in [0,1)$ such that 
\begin{equation*}
a = \frac{\alpha - \varepsilon_a}{N} 
\hspace{3ex}\text{and}\hspace{3ex}
b = \frac{\beta + \varepsilon_b}{N}.
\end{equation*}
(a) We may assume $\nu \ge 2$ and $N-1 \ge \beta-\alpha \ge \nu -2$, since otherwise (\ref{pre*}) holds trivially.
We first consider the case $\beta-\alpha \le N-2$. 
If $\nu$ points fall into $[a,b)$, then one of the three disjoint events that exactly $\nu$ points, $\nu-1$ points or $\nu-2$ points fall into $[\alpha/N, \beta/N)$ occurs.
Therefore 
\begin{equation*}
\begin{split}
 P_\nu :=& \PP \left( \bigwedge_{\ell=1}^\nu X_{\ell} \in [a,b) \right) = \frac{(\beta- \alpha)(\beta-\alpha -1) \cdots (\beta-\alpha -\nu +1)}{N(N-1) \cdots (N-\nu+1)} \\
 +& \nu \frac{(\beta- \alpha)(\beta-\alpha -1) \cdots (\beta-\alpha -\nu +2)}{N(N-1) \cdots (N-\nu+2)}  \frac{\varepsilon_a +
 \varepsilon_b}{N-\nu +1} \\ 
 +& \nu(\nu-1)\frac{(\beta- \alpha)(\beta-\alpha -1) \cdots (\beta-\alpha -\nu +3)}{N(N-1) \cdots (N-\nu+3)}  \frac{\varepsilon_a}{N-\nu +2} \frac{\varepsilon_b}{N-\nu +1}.
\end{split}
\end{equation*}
We have to verify that $P_\nu$ is at most $(b-a)^\nu = \big((\beta-\alpha +(\varepsilon_a+\varepsilon_b))/N \big)^\nu$. Since for fixed sum $\varepsilon_a+\varepsilon_b$ the product $\varepsilon_a\varepsilon_b$ is largest if $\varepsilon_a = \varepsilon_b =: \varepsilon$, we may confine ourselves to the latter case.
Put
\begin{equation*}
C:= \frac{(\beta- \alpha)(\beta-\alpha -1) \cdots (\beta-\alpha -\nu +3)}{N(N-1) \cdots (N-\nu+1)}
\end{equation*}
and define the function $f_\nu$ by
\begin{equation*}
\begin{split}
f_\nu(\varepsilon) &:= \\
&C \big[ (\beta - \alpha -\nu +2)(\beta - \alpha -\nu +1) + (\beta - \alpha - \nu + 2) 2\nu \varepsilon + 
\nu(\nu-1)\varepsilon^2 \big] \\
&\times \left( \frac{N}{\beta - \alpha + 2\varepsilon} \right)^\nu;
\end{split}
\end{equation*}
it suffices to show that $|f_\nu(\varepsilon)| \le 1$ for all $\varepsilon \in [0,1]$. It is easy to check that the zeros of the derivative 
$f_\nu'$ are 
$\varepsilon_1 = 1$ and $\varepsilon_2 \le 0$ if $\nu >2$ and $\varepsilon_1 = 1$ if $\nu =2$ . Hence $f_\nu$ takes its maximum in $[0,1]$ in $0$ or in $1$.
Now 
\begin{equation*}
f_\nu(0) = \frac{(\beta- \alpha)(\beta - \alpha -1)\cdots (\beta-\alpha -\nu +1)}{N(N-1) \cdots (N-\nu+1)}
 \left( \frac{N}{\beta - \alpha} \right)^\nu \le 1,
\end{equation*}
and from
\begin{equation*}
(\beta - \alpha -\nu +2)(\beta - \alpha -\nu +1) + (\beta - \alpha - \nu + 2) 2\nu + \nu(\nu-1) = (\beta - \alpha +2)(\beta - \alpha +1)
\end{equation*}
and $\beta-\alpha +2 \le N$ we obtain
\begin{equation*}
\begin{split}
f_\nu(1) &= \frac{(\beta-\alpha +2)(\beta-\alpha +1)(\beta- \alpha) \cdots (\beta-\alpha -\nu +3)}{N(N-1) \cdots (N-\nu+1)}
 \left( \frac{N}{\beta - \alpha + 2} \right)^\nu\\ 
 &\le 1.
 \end{split}
\end{equation*}
The remaining (simpler) case $\beta-\alpha = N-1$ can be verified similarly.

(b) Notice that we do not have to care about the  cases $a=0$ and $\nu=0$, since they are already covered by (a). 
Furthermore, due to (a) it suffices to show that 
\begin{equation*}
P^{(\sigma)}_C(k,\nu):= \PP \left(  \bigwedge^k_{\ell=\nu +1} X_\ell \in I^{(\sigma)}_2  \,\Big| \,  \bigwedge^{\nu}_{\ell=1} X_\ell \in I^{(\sigma)}_1  \right)
\le \delta \lambda(I_2^{(\sigma)})^{k-\nu}.
\end{equation*}

Let us first consider $\sigma = 0$: If $b\in G^1_N$, then $b= \beta/N$ and
\begin{equation*}
P^{(0)}_C(k,\nu) = \frac{\beta -\nu}{N-\nu}  \frac{\beta - (\nu+1)}{N-(\nu+1)} \cdots  \frac{\beta -(k-1)}{N-(k-1)}
\le \left( \frac{\beta}{N} \right)^{k-\nu} =  \big(\lambda(I^{(0)}_2) \big)^{k-\nu}.
\end{equation*}
If $b\notin G^1_N$, then $\beta \le N-1$ and
\begin{equation*}
\begin{split}
P^{(0)}_C(k,\nu) &\le \PP \left(  \bigwedge^k_{\ell=\nu +1} X_\ell \in \bigg[0, \frac{\beta+1}{N}\bigg)  \,\bigg| \,  \bigwedge^{\nu}_{\ell=1} X_\ell \in I^{(0)}_1  \right)\\
&= \frac{\beta +1-\nu}{N-\nu}  \frac{\beta +1 - (\nu+1)}{N-(\nu+1)} \cdots  \frac{\beta +1-(k-1)}{N-(k-1)}
\le \left( \frac{\beta}{N-1} \right)^{k-\nu}\\
&\le  \left( \frac{N}{N-1} \right)^{k-\nu} b^{k-\nu} \le \left( 1+\frac{1}{N-1} \right)^{N-1} b^{k-\nu}  
\le e \big(\lambda(I^{(0)}_2) \big)^{k-\nu}.
\end{split}
\end{equation*}

Let us now consider $\sigma = 1$: 
If $a,b \in G^1_N$, then we have $a=\alpha/N$, $b= \beta/N$, and it is easily verified that (\ref{*}) holds with $\delta =1$.

If $a \notin G^1_N$, then $\alpha \ge 1$ and 
\begin{equation*}
\begin{split}
&P^{(1)}_C(k,\nu) \le \PP \left(  \bigwedge^k_{\ell=\nu +1} X_\ell \in [0,a) \cup \bigg[\frac{\beta}{N}, 1 \bigg)  \,\bigg| \,  \bigwedge^{\nu}_{\ell=1} X_\ell \in I^{(1)}_1  \right)\\
&=  \frac{N-\beta + (\alpha - 1) - \nu}{N-\nu}  \frac{N-\beta + (\alpha -1) - (\nu + 1)}{N-(\nu+1)} \cdots \\ & \times \frac{N-\beta + (\alpha - 1) -(k-1)}{N-(k-1)}\\
&+ (k-\nu) \frac{N-\beta + (\alpha-1) -\nu}{N-\nu}  
\cdots  \frac{N-\beta + (\alpha-1) -(k-2)}{N-(k-2)} \frac{1-\varepsilon_a}{N-(k-1)}.\\
\end{split}
\end{equation*}
We want to prove that the last term is less or equal than 
\begin{equation*}
e \big( (N-\beta - \varepsilon_b + (\alpha -1) + (1- \varepsilon_a))/N \big)^{k-\nu}.
\end{equation*}
Obviously, it is enough to show this for the case $\varepsilon_b=1$.  
Put 
\begin{equation*}
C:= \frac{N-\beta + (\alpha-1) -\nu}{N-\nu}  
\cdots  \frac{N-\beta + (\alpha-1) -(k-2)}{N-(k-2)} \frac{1}{N-(k-1)}
\end{equation*}
and define the function $f_{k,\nu}$ by 
\begin{equation*}
f_{k,\nu}(\varepsilon) := C[N-\beta + (\alpha - 1) - (k-1) + (k-\nu)\varepsilon] \left( \frac{N}{N-\beta + \alpha - 2 + \varepsilon} \right)^{k-\nu};
\end{equation*}
it suffices to show $|f_{k,\nu}(\varepsilon)| \le e$ for all $\varepsilon \in [0,1]$. The only zero of $f_{k,\nu}'$ is at least $1$, hence $f_{k,\nu}$ takes its maximum in $[0,1]$ in $0$ or $1$. 
Now
\begin{equation*}
\begin{split}
f_{k,\nu}(0) =& \frac{N-\beta + (\alpha -1) -\nu}{N-\nu} \cdots \frac{N-\beta + (\alpha -1) - (k-1)}{N- (k-1)}\\
&\times  \left( \frac{N}{N-\beta + (\alpha -1) -1} \right)^{k-\nu}\\
\le& \left( \frac{N}{N-1} \right)^{k-\nu} \le e, 
\end{split}
\end{equation*}
and
\begin{equation*}
\begin{split}
f_{k,\nu}(1) =& \frac{N-\beta + (\alpha -1) -(\nu-1)}{N-\nu} \cdots \frac{N-\beta + (\alpha -1) - (k-2)}{N- (k-1)}\\ 
&\times \left( \frac{N}{N-\beta + (\alpha -1)} \right)^{k-\nu}\\
\le& \left( \frac{N}{N-1} \right)^{k-\nu} \le e. 
\end{split}
\end{equation*}
We now show that the constant $\delta$ in (\ref{*}) is optimal: Let $\sigma = 0$ and $\beta = N-1 = \alpha$, $\varepsilon_a=0$, 
$\varepsilon_b\in (0,1)$, $k=N$, and $\nu =1$. For this choice of parameters we get 
\begin{equation*}
\PP \left(  X_\ell \in I^{(0)}_1 \right) = b-a 
\hspace{3ex}\text{and}\hspace{3ex} 
P^{(0)}_C(N,1) = 1.
\end{equation*}
Therefore every $\delta$ that satisfies (\ref{*}) for every choice of $N$ and $\varepsilon_b$ has to fulfill
\begin{equation*}
\delta \ge \left( \frac{N}{N-1+\varepsilon_b} \right)^{N-1}  \ge  \left( 1+ \frac{1-2 \varepsilon_b}{N-1} \right)^{N-1};
\end{equation*}
since the last expression converges to $e^{1-2\varepsilon_b}$ for $N\to \infty$ and since we can choose $\varepsilon_b$ arbitrarily small, this implies $\delta \ge e$.

In the case $\sigma = 1$  we can consider a corresponding example: Let $\alpha= N-1 = \beta$, $\varepsilon_a=0$, $\varepsilon_b \in (0,1)$, $k=N$, and $\nu=1$. Again, it is easily verified that (\ref{*}) cannot hold for $\delta<e$ by choosing $N$ sufficiently large and $1-\varepsilon_b$ sufficiently small. 

This concludes the proof of Lemma \ref{one-d}.
\end{proof}

\begin{definition}\label{Def_neg_dep}
For $d\in \N$ we put 
$$\mathcal{C}^d_0 : = \{ [0,a) \,|\, a\in [0,1]^d\}
\hspace{3ex}\text{and}\hspace{3ex} 
\mathcal{D}^d_0 : = \{ B\setminus A \,|\, A,B \in \mathcal{C}^d_0 \}.$$
Let $\mathcal{S}\in \{\mathcal{C}^d_0, \mathcal{D}^d_0\}$. We say that 
the random points $X_1,\ldots,X_N$ in $[0,1)^d$ are \emph{$\mathcal{S}$-$\gamma$-negatively dependent} if for all $S\in \mathcal{S}$ the random variables 
\begin{equation*}
\1_S(X_1), \ldots, \1_S(X_N) \hspace{2ex}\text{are $\gamma$-negatively dependent.}
\end{equation*} 
\end{definition}

\begin{theorem}\label{Thm-d-LHS}
Let $d, d', d''\in \N_0$ such that $d= d' + d''$. 
Let $(X_n)_{n=1}^N$ be a LHS in $[0,1)^{d'}$ and $(Y_n)_{n=1}^N$ be independently randomized Monte Carlo points in $[0,1)^{d''}$. For $n=1, \ldots, N$ put $Z_n:= (X_n, Y_n)$. For $a, b \in [0,1)^d$ let $A:= [0,a)$, $B:= [0,b)$, and
$D:= B \setminus A$. 
Then the  random variables 
\begin{equation}\label{LHS_neg_dep}
\1_D(Z_1), \ldots, \1_D(Z_N) \hspace{2ex}\text{are $\gamma_{d'}$-negatively dependent,}
\end{equation}
where 
\begin{equation*}
\gamma_{d'} := \prod_{i=1}^{d'} \delta_i \hspace{2ex}\text{and} \hspace{2ex} \delta_i := \begin{cases} 1 & \mbox{if} \quad a_i, b_i \in G^1_N \quad \mbox{or} \quad a_i=0,\\ 
 e & \mbox{else}.\end{cases}
\end{equation*}
In particular, 
the random points $(Z_n)_{n=1}^N$ are $\mathcal{C}^d_0$-negatively as well as $\mathcal{D}^d_0$-$\gamma_{d'}$-negatively dependent.
\end{theorem}


\begin{proof}
Let $U=(u_n)^N_{n=1}$ be a family of uniformly distributed i.i.d. random points in $[0,1)^d$.
For $c\in \{0,1,\ldots, d\}$ we define the random point sets $\widehat{P}(c) = (\widehat{p}_n(c))^N_{n=1}$ by
\begin{equation*}
(\widehat{p}_n(c))_i = \begin{cases} \frac{\pi_i(n-1) + u_{n,i}}{N}& \mbox{for}\quad i\le c\\ 
\quad u_{n,i}& \mbox{for}\quad i>c,\end{cases}
\end{equation*} 
where $\pi_i$ is a randomly chosen permutation. Here the permutations $\pi_i$, $i\in [d]$, and the 
$u_n$, $n\in [N]$, are mutually independent. Notice that $\widehat{P}(c)$ is an MC point set for $c=0$, a $c$-dimensional
LHS padded by MC for $1\le c<d$, and a $d$-dimensional
LHS for $c=d$. Put 
\begin{equation*}
\gamma(c):= \prod^c_{i=1} \delta_i.
\end{equation*}
We first show via induction that for $c=0,1,\ldots,d$ the random variables 
\begin{equation}\label{star}
\1_D(\widehat{p}_1(c)), \ldots, \1_D(\widehat{p}_N(c))
\hspace{2ex}\text{are upper $\gamma(c)$-negatively dependent}.
\end{equation}
This is clearly satisfied for $c=0$, since the random variables 
$\1_D(\widehat{p}_1(0))$, $\ldots$, $\1_D(\widehat{p}_N(0))$ are even independent.
Now let $c\ge 1$ and assume that we have already showed (\ref{star}) for $c-1$. 
We use the shorthand $\widehat{P}:= \widehat{P}(c)$ and $\widetilde{P}:= \widehat{P}(c-1)$ and corresponding
notation for the random points in both sets. 
We denote by
$\mathcal{P}^*_c$ the orthogonal projection onto all coordinates except of the $c$th coordinate
(i.e., for $x\in \R^d$ we have $\mathcal{P}^*_c(x) = (x_1, \ldots, x_{c-1}, x_{c+1}, \ldots, x_d)$).
Note that $\mathcal{P}^*_c(\widehat{p}_j) = \mathcal{P}^*_c(\widetilde{p}_j)$ for all $j\in [N]$.
Furthermore, we put $A_c:= [0,a_c)$, $A^*_c := \mathcal{P}^*_c(A)$, $B_c:= [0,b_c)$, $B^*_c := \mathcal{P}^*_c(B)$,
and $D_c := B_c \setminus A_c$, $D^*_c := B^*_c \setminus A^*_c$. 
We have $\widehat{p}_j \in D$ if and only if one of the following two disjoint events occurs:
\begin{enumerate}
 \item $\mathcal{P}^*_c(\widehat{p}_j) = \mathcal{P}^*_c(\widetilde{p}_j) \in A^*_c$ 
 and $\widehat{p}_{j,c} \in D_c$,
\item $\mathcal{P}^*_c(\widehat{p}_j) = \mathcal{P}^*_c(\widetilde{p}_j) \in D^*_c$ 
 and $\widehat{p}_{j,c} \in B_c$.
\end{enumerate}
Since our random point distribution is symmetric, i.e., our random points are exchangeable, we get for  $k\in [N]$ 
\begin{equation}\label{condform}
\begin{split}
&\PP \left( \bigwedge^{k}_{j=1} \widehat{p}_j \in D \right)\\
&=\sum^k_{\nu = 0} {k \choose \nu} 
\PP \left( \bigwedge^{k}_{j=1} \widehat{p}_j \in D \wedge \bigwedge^{\nu}_{j=1} \mathcal{P}^*_c(\widetilde{p}_j ) 
\in A^*_c \wedge \bigwedge^{k}_{j=\nu + 1} \mathcal{P}^*_c(\widetilde{p}_j ) \in D^*_c \right)\\
&=\sum^k_{\nu=0} {k \choose \nu} \PP \left( \bigwedge^{\nu}_{j=1} \mathcal{P}^*_c(\widetilde{p}_j ) 
\in A^*_c \wedge \bigwedge^{k}_{j=\nu + 1} \mathcal{P}^*_c(\widetilde{p}_j ) \in D^*_c \right) \\
&\times \PP \left( \bigwedge^{\nu}_{j=1} \widehat{p}_{j,c} \in D_c \wedge \bigwedge^{k}_{j=\nu+1} \widehat{p}_{j,c} \in B_c  
\,\Big| \, \bigwedge^{\nu}_{j=1} \mathcal{P}^*_c(\widetilde{p}_j ) 
\in A^*_c \wedge \bigwedge^{k}_{j=\nu + 1} \mathcal{P}^*_c(\widetilde{p}_j ) \in D^*_c \right).
\end{split}
\end{equation}
Since different components of our random point set $\widehat{P}$ are 
mutually independent, for fixed $\nu\in \{0,1,\ldots,k\}$ the conditional probability in (\ref{condform}) is equal to 
\begin{equation*}
\PP \left( \bigwedge^{\nu}_{j=1} \widehat{p}_{j,c} \in D_c \wedge \bigwedge^{k}_{j=\nu+1} \widehat{p}_{j,c} \in B_c \right)
\end{equation*}
(provided the event on which we condition occurs with positive probability).
Because of this observation and since (\ref{condform}) also holds for $\widetilde{P}$ if we substitute  all occuring points $\widehat{p}_j$ by $\widetilde{p}_j$, we obtain
\begin{equation}\label{**}
\PP \left( \bigwedge_{j=1}^k \widehat{p}_j \in D \right) \le  \delta_c \PP \left( \bigwedge_{j=1}^k \widetilde{p}_j \in D \right)  
\le  \gamma(c) \lambda^d(D)^k,
\end{equation}
where the first inequality follows from
\begin{equation*}
\begin{split}
\PP \left( \bigwedge^{\nu}_{j=1} \widehat{p}_{j,c} \in D_c \wedge \bigwedge^{k}_{j=\nu+1} \widehat{p}_{j,c} \in B_c \right)
\le&  \delta_c (b_c - a_c)^\nu b^{k-\nu}_c\\
=& \delta_c  \PP \left( \bigwedge^{\nu}_{j=1} \widetilde{p}_{j,c} \in D_c \wedge \bigwedge^{k}_{j=\nu+1} \widetilde{p}_{j,c} \in B_c \right) 
\end{split}
\end{equation*}
which is valid due to (\ref{*}) for $\sigma=0$,  and the second inequality follows from our induction hypothesis. 

Now we show that $1_D(\widehat{p}_1(c)), \ldots, 1_D(\widehat{p}_N(c))$ are lower $\gamma(c)$-negatively
dependent. This holds if and only if
\begin{equation}\label{***}
 1_F(\widehat{p}_1(c)), \ldots, 1_F(\widehat{p}_N(c)) \hspace{2ex} \text{are upper $\gamma(c)$-negatively dependent,}
\end{equation}
where $F:= A \cup \big( [0,1)^d\setminus B \big)$.
We now verify (\ref{***}) by induction. Again, the statement is obvious for $c=0$. So let $c\ge 1$ and assume (\ref{***}) holds for $c-1$. As before we use
the notation $\widehat{P}$, $\widetilde{P}$, etc..
We have $\widehat{p}_j \in F$ if and only if one of the following three disjoint events occurs:
\begin{enumerate}
 \item $\mathcal{P}^*_c(\widehat{p}_j) = \mathcal{P}^*_c(\widetilde{p}_j) \in A^*_c$ 
 and $\widehat{p}_{j,c} \in A_c \cup [b_c,1)$,
\item $\mathcal{P}^*_c(\widehat{p}_j) = \mathcal{P}^*_c(\widetilde{p}_j) \in D^*_c$ 
 and $\widehat{p}_{j,c} \in [b_c,1)$,
\item $\mathcal{P}^*_c(\widehat{p}_j) = \mathcal{P}^*_c(\widetilde{p}_j) \in [0,1)^{d-1} \setminus B^*_c$ 
 and $\widehat{p}_{j,c} \in [0,1)$.
\end{enumerate}
Using the standard notation for multinomial coefficients, we obtain for $k\in [N]$
\begin{equation}\label{2.11'}
 \begin{split}
  &\PP \left( \bigwedge^k_{j=1} \widehat{p}_j \in F \right) = \sum_{0 \le \nu_1 \le \nu_2 \le k} {k \choose \nu_1, \nu_2 - \nu_1, k-\nu_2 } \\
 &\times \PP \left( \bigwedge^{\nu_1}_{j=1} \mathcal{P}^*_c(\widetilde{p}_j) \in A^*_c \wedge \bigwedge^{\nu_2}_{j=\nu_1+1} 
 \mathcal{P}^*_c(\widetilde{p}_j) \in D^*_c \wedge \bigwedge^{k}_{j=\nu_2+1} \mathcal{P}^*_c(\widetilde{p}_j) \in [0,1)^d\setminus B^*_c
 \right)\\
 &\times \PP \bigg( \bigwedge^{\nu_1}_{j=1} \widehat{p}_{j,c} \in A_c \cup [b_c,1) \wedge \bigwedge^{\nu_2}_{j=\nu_1+1} 
 \widehat{p}_{j,c} \in [b_c,1) \wedge \bigwedge^{k}_{j=\nu_2+1} \widehat{p}_{j,c} \in [0,1) \, \bigg|\\
 &\,\bigg|\, \bigwedge^{\nu_1}_{j=1} \mathcal{P}^*_c(\widetilde{p}_j) \in A^*_c \wedge \bigwedge^{\nu_2}_{j=\nu_1+1} 
 \mathcal{P}^*_c(\widetilde{p}_j) \in D^*_c \wedge \bigwedge^{k}_{j=\nu_2+1} \mathcal{P}^*_c(\widetilde{p}_j) \in [0,1)^d\setminus B^*_c
 \bigg).
 \end{split}
\end{equation}
For fixed $\nu_1, \nu_2$ the conditional probability appearing in the sum in (\ref{2.11'}) is equal to 
\begin{equation*}
 \PP \left( \bigwedge^{\nu_1}_{j=1} \widehat{p}_{j,c} \in A_c \cup [b_c,1) \wedge \bigwedge^{\nu_2}_{j=\nu_1+1} 
 \widehat{p}_{j,c} \in [b_c,1) \right)
\end{equation*}
(provided that the event on which we condition occurs with positive probability).
Since this observation and (\ref{2.11'}) hold also for $\widetilde{P}$ if we substitute all occuring points $\widehat{p}_j$ by $\widetilde{p}_j$,
the inequality
\begin{equation}\label{****}
 \PP \left( \bigwedge^k_{j=1} \widehat{p}_j \in F \right) \le \delta_c \PP \left( \bigwedge^k_{j=1} \widetilde{p}_j \in F \right) 
 \le \gamma(c) \lambda^d(F)^k
\end{equation}
follows from 
\begin{equation*}
\begin{split}
 \PP \bigg( \bigwedge^{\nu_1}_{j=1} \widehat{p}_{j,c} \in A_c \cup [b_c,1) &\wedge \bigwedge^{\nu_2}_{j=\nu_1+1} 
\widehat{p}_{j,c} \in [b_c,1) \bigg) \le \delta_c (a_c + (1-b_c))^{\nu_1} (1-b_c)^{\nu_2}\\
 &= \delta_c \PP \left( \bigwedge^{\nu_1}_{j=1} \widetilde{p}_{j,c} \in A_c \cup [b_c,1) \wedge 
 \bigwedge^{\nu_2}_{j=\nu_1+1} \widetilde{p}_{j,c} \in [b_c,1) \right),
\end{split}
\end{equation*}
which holds true due to (\ref{*}) for $\sigma=1$, and our induction hypothesis.

This concludes the proof of Theorem  \ref{Thm-d-LHS}.
\end{proof}

Studying the proof above  it is easy to see that the following generalization of Theorem~\ref{Thm-d-LHS} is valid.

\begin{theorem}\label{theo}
Let $a, b \in [0,1)^d$ and put $A:= [0,a)$, $B:= [0,b)$, and $D:=B \setminus A$. 
Let $(Z_n)_{n=1}^N$ be a set of (not necessarily independent) random points in $[0,1)^d$ that satisfies the following two conditions:
\begin{itemize}
\item[{\rm (i)}] Different components of the random point set are mutually independent.
\item[{\rm (ii)}] For all $i\in [d]$ and  for
$I^{(0)}_{1,i} := [a_i,b_i)$, $I^{(0)}_{2,i} := [0,b_i)$ and $I^{(1)}_{1,i} := [b_i,1)$, $I^{(1)}_{2,i} := [0,a_i) \cup [b_i,1)$ there exists a $\delta_i >0$ such that for all $\sigma\in \{0,1\}$, $k\in \{0,1, \ldots, N\}$, $\nu \in \{ 0,1, \ldots, k\}$, and all $J\subseteq [N]$,  
$J_\nu$, $J_{k-\nu}\subseteq J$ with 
$|J|=k$, $|J_\nu|= \nu$, $|J_{k-\nu}| = k-\nu$ and $J_\nu \cap J_{k-\nu} = \emptyset$,
one has
\begin{equation}\label{as*}
 \PP \left( \bigwedge_{\ell \in J_\nu} Z_{\ell,i} \in I^{(\sigma)}_{1,i} \wedge  \bigwedge_{\ell\in J_{k-\nu}} Z_{\ell,i} \in I^{(\sigma)}_{2,i} \right) 
 \le \delta_i \lambda( I^{(\sigma)}_{1,i})^{\nu} \lambda( I^{(\sigma)}_{2,i})^{k-\nu}.
\end{equation}
\end{itemize}
Then the  random variables 
\begin{equation}\label{A_neg_dep}
\1_D(Z_1), \ldots, \1_D(Z_N) \hspace{2ex}\text{are $\gamma_d$-negatively dependent, where} \hspace{2ex} \gamma_d := \prod_{i=1}^d \delta_i.
\end{equation}
\end{theorem}

\section{Probabilistic Discrepancy Bounds}
\label{PDB}

Now we consider the star discrepancy $D_N^*(X)$ (as defined in the introduction) of $\mathcal{D}^d_0$-$\gamma$-negatively dependent random points $X=(X_n)_{n=1}^N$ (cf. Definition  \ref{Def_neg_dep}).

To ``discretize'' the star discrepancy, we define $\delta$--covers 
as in~\cite{DGS05}: 

\begin{definition}
For any $\delta\in(0,1]$ a finite set $\Gamma$ of points in $[0,1)^d$ is called a \emph{$\delta$--cover} of $[0,1)^d$, if for every $y\in [0,1)^d$ there exist $x,z\in\Gamma\cup\{0\}$ such that $x\leq y\leq z$ and $\lambda^d([0,z])-\lambda^d([0,x])\leq\delta$. The number $\mathcal{N}(d,\delta)$ denotes the smallest cardinality of a 
$\delta$--cover of $[0,1)^d$. 
\end{definition}

The following theorem was stated and proved in~\cite{Gne08a}.

\begin{theorem}\cite[Thm.1.15]{Gne08a}
\label{bracketing}
For any $d\geq 1$ and $\delta\in(0,1]$ we have
\begin{equation*}
\mathcal{N}(d,\delta)\leq 2^d \frac{d^d}{d!}(\delta^{-1}+1)^d.
\end{equation*}
\end{theorem}

Notice that due to Stirling's formula we have $d^d/d! \le e^d/\sqrt{2\pi d}$. 
Furthermore, it is easy to verify that in the case $d=1$ the identity 
\begin{equation}\label{cover_d=1}
\mathcal{N}(1,\delta) = \lceil \delta^{-1} \rceil
\end{equation}
is established with the help of the $\delta$-Cover $\Gamma := \{1/ \lceil \delta^{-1} \rceil, 2/ \lceil \delta^{-1} \rceil, \ldots,1\}$.

With the help of $\delta$-covers the star discrepancy can be approximated in the following sense.

\begin{lemma}\label{delta_approx}
Let $P\subset [0,1)^d$ be an $N$-point set, $\delta >0$, and $\Gamma$ be a $\delta$-cover  of $[0,1)^d$. Then
\begin{equation*}
D^*_N(P) \le \max_{x\in \Gamma} D_N(P, [0,x)) + \delta.
\end{equation*}
\end{lemma}

The proof of Lemma \ref{delta_approx} is straightforward, cf., e.g., \cite[Lemma~3.1]{DGS05}.

We are ready to state and prove our main result, a general probabilistic discrepancy bound for 
$\mathcal{D}^d_0$-$\gamma$-negatively dependent random points.

\begin{theorem}\label{main_theo}
Let $d, N\in \N$ and $\rho \in [0,\infty)$. Let $X= (X_n)_{n\in[N]}$ be a set of $\mathcal{D}^d_0$-$e^{\rho d}$-negatively dependent random points in $[0,1)^d$  
such that each $X_n$ is uniformly distributed. 
Then for every $c>0$
\begin{equation}\label{disc_bound_GH}
D^*_N(X) \leq c\sqrt{\frac{d}{N}}
\end{equation}
holds with probability at least $1-e^{-(1.6741\cdot c^2 -10.7042- \rho)\cdot d}$, implying that for every $q\in (0,1)$ 
\begin{equation}\label{disc_bound_AH+}
D^*_N(X) \leq 0.7729 \sqrt{ 10.7042 + \rho + \frac{\ln \big( (1-q)^{-1} \big) }{d}} \sqrt{\frac{d}{N}}
\end{equation}
holds with probability at least $q$. 
\end{theorem}

Notice that in the special case of Monte Carlo points $X_1,\ldots,X_n$ in Theorem~\ref{main_theo}
we have $\rho =0$ and our bound \eqref{disc_bound_AH+} improves on the main result of \cite{AH14}, that is, Theorem~1 in that paper.
This is further illustrated by the next corollary. It follows immediately from Theorem~\ref{main_theo}, since we may take $\rho=0$ for Monte Carlo points and $\rho=1$ for Latin hypercube samples (padded by Monte Carlo or not), see Theorem~\ref{Thm-d-LHS}. Notice in particular the quantitative improvements in the constant in \eqref{aist*} compared to the constant $9.65$ provided in \cite{Ais11} and in the discrepancy bounds \eqref{aisthof*} compared to the explicit discrepancy bounds for Monte Carlo point sets presented in \cite{AH14} in the table on p.~1374.

\begin{corollary}\label{Cor_Prob_Est}
Let $d, N\in \N$ and let $X= (X_n)_{n\in[N]}$ be a random point set in $[0,1)^d$.
\begin{enumerate}
\item If $X$ is a Monte Carlo point set, then there exists a realization $P\subset [0,1)^d$ of $X$ such that
\begin{equation}\label{aist*}
D^*_N(P) \leq 2.5287 \cdot\sqrt{\frac{d}{N}}.
\end{equation}
The probability that $X= (X_n)_{n\in[N]}$ satisfies
\begin{equation}\label{aisthof*}
D^*_N(X) \leq 3\cdot\sqrt{\frac{d}{N}}
\hspace{3ex}\text{and}\hspace{3ex}
D^*_N(X) \leq 4\cdot\sqrt{\frac{d}{N}}
\end{equation}
is at least $0.987256$ and $0.999999$, respectively. 
\item If $X$ is a Latin hypercube sample or a Latin hypercube sample padded by Monte Carlo, then there exists a realization $P\subset [0,1)^d$ of $X$ such that
\begin{equation}
D^*_N(P) \leq 2.6442 \cdot\sqrt{\frac{d}{N}}.
\end{equation}
The probability that $X= (X_n)_{n\in[N]}$ satisfies
\begin{equation}
D^*_N(X) \leq 3\cdot\sqrt{\frac{d}{N}}
\hspace{3ex}\text{and}\hspace{3ex}
D^*_N(X) \leq 4\cdot\sqrt{\frac{d}{N}}
\end{equation}
is at least $0.965358$ and $0.999999$, respectively. 
\end{enumerate}
Estimate \eqref{aist*} implies
\begin{equation}
N^*(\varepsilon ,d) \le \left\lceil 6.3944 \cdot d \varepsilon^{-2} \right\rceil
\hspace{3ex}\text{for all $\varepsilon >0$.}
\end{equation}
\end{corollary}

We may use Theorem~\ref{main_theo} to prove similar corollaries for other random samples than Monte Carlo
point sets or (padded) Latin hypercube samples; cf. also \cite{WGH19}.

\begin{remark}\label{Rem_DDG18}
We already mentioned in the introduction that the probabilistic discrepancy bound \eqref{disc_bound_GH}
is sharp for Monte Carlo point sets. As shown in \cite{DDG18}, the same is the case for Latin hypercube samples: 
 There exists a constant $K>0$ such that for $d\ge 2$ and $N \ge 1600 d$ the expected star discrepancy of a Latin hypercube sample $X$ is bounded from below by
\begin{equation}\label{ddg1}
\E [D^*_N(X)] \ge K \sqrt{\frac{d}{N}}
\end{equation}
and additionally we have the probabilistic discrepancy bound
\begin{equation}\label{ddg2}
\PP \! \left( D^*_N(X) < K\sqrt{\frac{d}{N}} \right) \le \exp(-\Omega(d)),
\end{equation}
see \cite[Theorem~2]{DDG18}.
Nevertheless, recall that Latin hypercube samples have a big advantage over Monte Carlo samples, namely that their one-dimensional projections are more evenly distributed, cf. Remark~\ref{1_dim_projections}. 
\end{remark}


\begin{proof}[Proof of Theorem \ref{main_theo}]
We adapt the line of proof of \cite[Theorem~1]{Ais11} and employ the bounds on the size of minimal $\delta$-covers from Theorem~\ref{bracketing}, dyadic chaining and large deviation bounds of Hoeffding- and Bernstein-type (but this time the ones for sums of $\gamma$-negatively dependent random variables in Theorem \ref{Hoeffding} and \ref{Bernstein}).
For $a,b\in [0,1]^d$ with $a\le b$ we write
\begin{equation*}
\Delta(a, b) := [0,b)\setminus [0,a).
\end{equation*}
We start by putting $\mu:= 13$ and 
\begin{equation*}
c_\mu := \sum^\infty_{\ell =0} \left( \sqrt{\frac{\mu+1}{2\mu}} \right)^\ell
= \frac{1}{1- \sqrt{\frac{\mu + 1}{2\mu}}}.
\end{equation*}
Let $c_0>0$ be given, and let $c_1>0$ be a constant; the value of $c_1$ will be determined later in the proof. 
Let $K$  be the smallest natural number that satisfies $K\ge \mu$ and 
\begin{equation}\label{def_K}
\frac{1}{\sqrt{K 2^{K}}} \le c_0 c_1 c_{\mu} \sqrt{\frac{d}{N}}.
\end{equation}
We choose for each $\mu \le k \le K$ a $2^{-k}$-cover $\Gamma_k$ of minimum size.
Furthermore, we put $\Gamma_{\mu-1} := \{0\}$.

Let $P$ be an arbitrary realization of $X$. Due to Lemma~\ref{delta_approx} we can choose for each test box $[0,y) \subseteq [0,1)^d$ an 
$a_K \in \Gamma_K \cup \{0\}$ such that 
\begin{equation*}
D_N(P, [0,y)) \le D_N(P, [0,a_K)) + 2^{-K}.
\end{equation*}
If $K>\mu$, we additionally choose for $k= K-1, \ldots, \mu$  points $a_k = a_k(a_{k+1}) \in \Gamma_k \cup \{0\}$
recursively, depending only on the previously chosen point $a_{k+1}$, such that $a_k \le a_{k+1}$ and  
\begin{equation}\label{lambda_Delta}
\lambda^d(\Delta(a_k, a_{k+1})) \le 2^{-k}.
\end{equation} 
Finally, we put $a_{\mu-1} = a_{\mu-1}(a_\mu) := 0$ and get $\Delta(a_{\mu-1},a_\mu) = [0,a_\mu)$.
Notice that $[0,a_K) = \cup_{k=\mu}^K \Delta(a_{k-1}, a_k)$.  
Hence we have 
\begin{equation}\label{disc_est}
D_N(P,[0,y)) \le \sum_{k=\mu}^K D_N(P,\Delta(a_{k-1}, a_k)) + 2^{-K}.
\end{equation}
Let us now for $\mu \le k\le K$ define the sets $\mathcal{A}_k$ by
\begin{equation*}
\mathcal{A}_k := \{ \Delta(a_{k-1}(b),b)  \,|\,  b \in \Gamma_k\}.
\end{equation*}
Clearly, $|\mathcal{A}_k| \le |\Gamma_k|$.
Moreover, we define events $E_k$ by 
\begin{equation*}
E_\mu := \left\{ \max_{\Delta_{\mu} \in \mathcal{A}_\mu} D_N(X, \Delta_\mu) \le c_0 \sqrt{\frac{d}{N}} \right\}
\end{equation*}
for $k=\mu$ and 
\begin{equation*}
E_k := \left\{  \max_{\Delta_{k} \in \mathcal{A}_k} D_N(X, \Delta_k) \le c_0 c_1 \sqrt{\frac{k-1}{2^{k-1}}}
\sqrt{\frac{d}{N}} \right\}
\end{equation*}
for $\mu +1 \le k \le K$. We put
\begin{equation*}
E:= \bigcap^K_{k=\mu} E_k.
\end{equation*}
Let us assume that the event $E$ has occured. Then, due to (\ref{disc_est}) and (\ref{def_K}), the realization $P$ of $X$ satisfies for an arbitrary test box $[0,y) \subseteq [0,1)^d$
\begin{equation*}
\begin{split}
&D_N(P, [0,y))\\
\le &c_0 \left( 1 + c_1 \sum_{k=\mu + 1}^K  \sqrt{\frac{k-1}{2^{k-1}}} \,\right) \sqrt{\frac{d}{N}} + 2^{-K}\\ 
\le &c_0  \left( 1 + c_1 \left( \sqrt{\frac{\mu}{2^\mu}}\,
\sum_{j=0}^{K-\mu-1}  \sqrt{\frac{\mu+j}{2^j\mu}}  +  c_\mu \sqrt{\frac{K}{2^K}} \,\right) \right) 
\sqrt{\frac{d}{N}}\\
\le &c_0  \left( 1 + c_1 \sqrt{\frac{\mu}{2^\mu}}
 \left(\sum_{j=0}^{K-\mu-1}  \left( \sqrt{\frac{\mu+1}{2\mu}} \, \right)^j +   \sqrt{\frac{\mu + (K-\mu)}{\mu 2^{K-\mu}}} c_\mu\,\right) \right) 
\sqrt{\frac{d}{N}}\\
\le &c_0  \left( 1 + c_1 \sqrt{\frac{\mu}{2^\mu}}\,
 \left( \sum_{j=0}^{K-\mu-1}  \left( \sqrt{\frac{\mu+1}{2\mu}} \, \right)^j +   \left( \sqrt{\frac{\mu + 1}{2 \mu}} \,\right)^{K-\mu} c_\mu\,\right) \right) 
\sqrt{\frac{d}{N}}\\
= &c_0 \left( 1 + c_1 c_\mu \sqrt{\frac{\mu}{2^\mu}} \,\right) \sqrt{\frac{d}{N}}.
\end{split}
\end{equation*}

Let us now derive a lower bound for the probability $\PP(E)$. 
For $k=\mu$ we may use Theorem \ref{Hoeffding} and a simple union bound to obtain
\begin{equation}\label{hoffi}
\PP(E^c_\mu) \le 2|\Gamma_\mu| e^{\rho d} e^{-2c_0^2 d}. 
\end{equation}
For $\mu + 1\le k \le K$ we first use the definition of $K$, cf. (\ref{def_K}), to get  the estimate
\begin{equation*}
\sqrt{(k-1)2^{1-k}} \le \sqrt{(K-1)2^{K-1}}\, 2^{1-k} < \left( c_0c_1c_\mu \sqrt{ \frac{d}{N} } \right)^{-1} 2^{1-k}.
\end{equation*}
Theorem \ref{Bernstein} together with (\ref{lambda_Delta}) and the latter estimate gives us 
\begin{equation}\label{bernie}
\PP(E^c_k) \le 2 |\Gamma_k| e^{\rho d} \exp \left( - \frac{c_0^2c_1^2 (k-1) d}{2(1+1/3c_\mu)} \right).
\end{equation}
We have 
\begin{equation*}
\PP(E) = 1 - \PP(E^c) \ge 1 - \sum_{k=\mu}^K \PP(E^c_k) . 
\end{equation*}
We put
\begin{equation*}
\tau_\mu := \frac{c^2_1}{4(1+1/3c_\mu)}
\end{equation*}
and $\sigma:= \mu - \ln(2(2^\mu +1)) -1$. 
From now on we assume that
\begin{equation*}
c_0 \ge \sqrt{(\mu+\rho -\sigma)/2}.
\end{equation*}

Let us first consider the case $d=1$. Due to (\ref{cover_d=1}) we have
$|\Gamma_k| = 2^{k}$ for $k=\mu,\ldots,K$. Hence we get from (\ref{hoffi}) and (\ref{bernie}) for $k=\mu$
\begin{equation*}
\PP(E^c_\mu) \le \exp \left( -2c_0^2 + \rho + (\mu+1)\ln(2) \right) 
\end{equation*}
and for $\mu+1 \le k \le K$
\begin{equation*}
\PP(E^c_k) \le \exp \left( - (2c_0^2\tau_\mu -\ln(2))(k-1) + \rho + 2\ln(2) \right).
\end{equation*}
Hence we obtain 
\begin{equation*}
\begin{split}
&\PP(E) 
\ge 1 -  \left(  e^{-(2c_0^2 -\rho - (\mu+1)\ln(2))} + e^{\rho +2\ln(2)}
\sum_{k=\mu+1}^{K}e^{ - \left(2c_0^2\tau_\mu - \ln(2) \right)(k-1)} \right)\\
= &1 -  \left(  e^{-(2c_0^2 -\rho - (\mu+1)\ln(2))} + e^{\rho +2\ln(2)} e^{ - \left(2c_0^2\tau_\mu - \ln(2) \right)\mu}
\sum_{j=0}^{K-\mu-1}e^{ - \left(2c_0^2\tau_\mu - \ln(2) \right)j} \right)\\
= &1 -   e^{-(2c_0^2 -\rho - (\mu+1)\ln(2))} \left( 1 + \frac{e^{ - 2c_0^2(\mu \tau_\mu - 1) + \ln(2)}}{1 - e^{ - \left(2c_0^2\tau_\mu - \ln(2) \right)}} \right)\\
\ge &1 -   e^{-(2c_0^2 -\rho - (\mu+1)\ln(2))} \left( 1 + \frac{e^{ - (\mu+\rho -\sigma)(\mu \tau_\mu - 1) + \ln(2)}}{1 - e^{ - (\mu+\rho-\sigma)\tau_\mu + \ln(2) }} \right).
\end{split}
\end{equation*}
Choosing $\tau_\mu = 0,0887$, 
we get for any $\rho \ge 0$
\begin{equation*}
1 + \frac{e^{ - (\mu+\rho -\sigma)(\mu \tau_\mu - 1) + \ln(2)}}{1 - e^{ - (\mu+\rho-\sigma)\tau_\mu + \ln(2) }} < e.
\end{equation*}
Since
\begin{equation*}
e^{-(2c_0^2 -\rho - (\mu+1)\ln(2))} 
< \frac{1}{e} e^{-(2c_0^2 -\rho - \mu+\sigma)} ,
\end{equation*}
the choice $c_0 = \sqrt{(\mu+\rho -\sigma)/2}$ leads to $\PP(E)>0$.

We now consider the case $d\ge 2$. Due to Theorem \ref{bracketing} we have
\begin{equation*}
|\Gamma_k| \le \frac{1}{\sqrt{2\pi d}} (2e)^d (2^k + 1)^d \hspace{3ex}\text{for $k=\mu,\ldots,K$.}
\end{equation*}
Hence we get from (\ref{hoffi}) and (\ref{bernie}) for $k=\mu$
\begin{equation*}
\PP(E^c_\mu) \le \sqrt{\frac{2}{\pi d}} \,e^{-(2c_0^2 - \mu - \rho + \sigma)d} 
\end{equation*}
and for $\mu+1 \le k \le K$
\begin{equation*}
\begin{split}
\PP(E^c_k) &\le \sqrt{\frac{2}{\pi d}}\, (2e)^d 2^{kd} (1+ 2^{-k})^d e^{\rho d} \exp \left( - \frac{c_0^2c_1^2(k-1)d}{2(1+1/3c_\mu)} \right)\\
&\le \sqrt{\frac{2}{\pi d}}\, e^{\left( 1+2\ln(2) + \vartheta + \rho \right)d}\exp \left( - (2c_0^2\tau_\mu -\ln(2))(k-1) d\right),
\end{split}
\end{equation*}
where $\vartheta:= \ln(1+2^{-\mu-1})$. Put $\zeta:= 2\ln(2) + \vartheta$. 
Then we obtain 
\begin{equation*}
\begin{split}
\PP(E) &\ge 1 -  \sqrt{\frac{2}{\pi d}}  \bigg(  e^{-(2c_0^2 -\mu -\rho + \sigma)d} \\
&\quad+ e^{(1 + \zeta +\rho)d} e^{ - \left(2c_0^2\tau_\mu - \ln(2) \right)\mu d}
\sum_{j=0}^{K-\mu-1}e^{ - \left(2c_0^2\tau_\mu - \ln(2) \right)jd} \bigg)\\
\ge &1 -   \sqrt{\frac{2}{\pi d}}\,e^{-(2c_0^2 -\mu -\rho +\sigma)d} \left( 1 + \frac{e^{ - (2c^2_0(\mu \tau_\mu - 1) + (1-\ln(2))\mu -1 -\zeta -\sigma)d}}{1 - e^{ - (2c^2_0\tau_\mu - \ln(2))d }} \right)\\
\ge &1 -   \sqrt{\frac{2}{\pi d}}\,e^{-(2c_0^2 -\mu -\rho +\sigma)d} \left( 1 + \frac{e^{ - ( (\mu +\rho -\sigma)(\mu \tau_\mu - 1) + (1-\ln(2))\mu -1 -\zeta -\sigma )d}}{1 - e^{ - ( (\mu +\rho -\sigma)\tau_\mu - \ln(2) )d }} \right).
\end{split}
\end{equation*}
Choosing $\tau_\mu = 0,0887$ 
as in the case $d=1$, we get 
\begin{equation}\label{quetsch_d}
1 + \frac{e^{ - ( (\mu +\rho -\sigma)(\mu \tau_\mu - 1) + (1-\ln(2))\mu -1 -\zeta -\sigma )d}}{1 - e^{ - ( (\mu +\rho -\sigma)\tau_\mu - \ln(2) )d }} < \sqrt{\frac{\pi d}{2}};
\end{equation}
this can easily be checked for $d= 2$ and -- since the left hand side of inequality (\ref{quetsch_d}) is monotonic decreasing in $d$, while the right hand side is monotonic increasing -- holds therefore for all $d\ge 2$.
Hence the choice $c_0 = \sqrt{(\mu+\rho -\sigma)/2}$ leads to $\PP(E)>0$.
\end{proof}

\subsection*{Acknowledgment}
The authors thank Marcin Wnuk and two anonymous referees for valuable comments.

Part of the work of Michael Gnewuch was done while he was a research fellow and a visitor at the
 School of Mathematics and Statistics of
the University of New South Wales in Sydney and a ``Chercheur Invit\'e'' at the
Laboratoire d'Informatique (LIX) of
\'Ecole Polytechnique in Paris, France. He acknowledges support from the
Australian Research Council ARC and thanks his hosts Josef Dick, Frances Y. Kuo, Ian H. Sloan, and Benjamin Doerr for their hospitality.

Another part of his work was done while
he visited special semesters and programs at the Radon Institute for Computational and Applied Mathematics (RICAM) in Linz, Austria, the Institute of Computational and Experimental Mathematics (ICERM) of Brown University in Providence, USA, and the Erwin Schr\"odinger Institute (ESI) in Vienna, Austria.



\end{document}